\documentstyle[12pt,amstex,amsfonts]{amsart}

\input epsf

\newtheorem{conjecture}{Conjecture}[section]
\newtheorem{theorem}{Theorem}[section]

\newtheorem{definition}{Definition}[section]

\newtheorem{lemma}{Lemma}[section]

\newenvironment{proof}{\noindent {\bf {Proof:}}}{$\Box$}

\numberwithin{equation}{section}

\begin{document} 

\title[Link invariants and the non-invertibility of links]
{Finite type link invariants and \\
the non-invertibility of links} 
\author{Xiao-Song Lin}
\address{Department of Mathematics, University of California, Riverside,
CA 92521}
\email{xl@@math.ucr.edu}
\date{First version: January 22, 1996. Revised: February 25, 1996.}
\thanks{Research supported in part by NSF}

\begin{abstract} 
We show that a variation of Milnor's $\bar\mu$-invariants, the so-called 
Campbell-Hausdorff invariants introduced recently by Stefan 
Papadima, are of finite type with respect to {\it marked singular links}.
These link invariants are stronger than quantum invariants in the sense that
they detect easily the non-invertibility of links with more than one 
components. 
It is still open whether some effectively computable knot invariants,
e.g.~ finite type knot invariants (Vassiliev 
invariants), could detect the non-invertibility of knots.    
\end{abstract}

\maketitle

\section{Introduction} 

It is well-known that the indeterminacy of Milnor's $\bar\mu$-invariants
\cite{Mi1,Mi2} is inadequate for various link classification problems (see, 
for example, \cite{Le,HaLi}). Moreover, 
the indeterminacy of $\bar\mu$-invariants
depends on individual 
links so that
a fixed $\bar\mu$-invariant takes values in different cyclic groups on
different links. This awkwardness is the reason why in 
\cite{Ba2,Li}, where one studies whether $\bar\mu$-invariants are 
of finite type, one has to drop off the 
indeterminacy at all to consider only string links. 

Recently, based on the work of M. Markl and S. Papadima \cite{MaPa} on an
algebraic parameterization for the set of (Mal'cev completed) fundamental
groups of the spaces with fixed first two Betti numbers ({\it a moduli space
of $h_{1,2}$-marked groups}), S. Papadima introduced a new class of link 
concordance invariants \cite{Pa}. For links with fixed number of components,
these so-called Campbell-Hausdorff invariants take 
values in a fixed moduli space. These invariants can be thought of as 
refinements of 
$\bar\mu$-invariants in that the latter are derived using the Magnus expansion
of free groups while the former using the Campbell-Hausdorff expansion (and 
hence the name). 

The moduli space of $h_{1,2}$-marked groups is a tower of 
orbit spaces of vector spaces of derivations of certain free nilpotent 
Lie algebras under the action
of a prounipotent group. 
To obtain the Campbell-Hausdorff invariants, one applies the 
Campbell-Hausdorff expansion to the quotients of the fundamental group
of a link complement by its lower central series. The indeterminacy is
thus reflected in the action of that prounipotent group. In 
\cite{Li}, we have argued that, regardless of the indeterminacy, link 
concordance invariants obtained from {\it any} expansion of free groups are of
finite type. The difficult now is that the moduli space where the 
Campbell-Hausdorff invariants live is not a vector space, although some
linear structures do survive there. To overcome 
this difficulty, we introduce the notion of {\it marked singular links}, which
are the usual singular links with a mark on each component away from double 
points. Under isotopy of marked singular links, the marks are only 
allow to 
move without passing through double points. Of course, isotopy classes of 
marked links are in one-one correspondence with isotopy classes of unmarked
links. The Campbell-Hausdorff link invariants can be used to derive
invariants of marked singular links. For a fixed Campbell-Hausdorff invariant,
we show that the derived marked singular link invariant vanishes when the 
number of double points is large enough. In that sense, we say that
the Campbell-Hausdorff invariants are of finite type with respect to 
marked singular links. It will certainly be interesting to explore 
systematically the 
combinatorics of link invariants of finite type with respect to marked
singular links. 

Being of finite type (with respect to unmarked singular links) is a 
property of 
link invariants obeyed by the 
coefficients of appropriate power series expansions of all link invariants 
derived from finite dimensional representations of simple Lie algebras
through the theory of quantum groups 
(generalizations of the Jones polynomial) 
\cite{Ba1,Bi,BiLi}. Since then, there has been a question about whether 
(linear combinations of) those finite type link invariants derived from 
quantum groups (quantum invariants) are the only finite type link invariants.
Notice that quantum invariants have an additional property that their values
will not change when the orientation of each component of a given link is
reversed due to a certain kind of self-duality of finite dimensional 
representations of simple Lie algebras. So one was asked whether 
finite type link invariants can detect the non-invertibility of a link. 
Recently, P. Vogel was able to answer the former question in negative 
through a very complicated combinatorial manipulation of chord diagrams. 
On the 
other hand, it is very simple to come up with examples of Campbell-Hausdorff
invariants which detect the non-invertibility of a link with more than one 
components. Therefore, such a Campbell-Hausdorff invariant can not be 
determined by quantum invariants on links with more than one 
components. It is still not known whether finite type knot invariants can
detect the non-invertibility of knots (we conjectured that the answer is no 
several years ago). 
 
Notice that the Campbell-Hausdorff invariants are derived directly from the 
fundamental
groups of links complements, while it has been a question whether the quantum
link invariants have anything to do with the fundamental groups of
link complements. We have shown here some of the similarities and differences 
of the Campbell-Hausdorff invariants and the quantum invariants. This,
we think, is of some significance.  

We will review the background in Section 2. The statement of the 
main result and its proof are in Section 3. In Section 4, we include a
proof of the well-known fact that quantum invariants can not detect the 
non-invertibility of a link. Section 5 contains an example
and some comments.

We are very grateful to Dror Bar-Natan for a critical reading of the original 
version of the paper. His comments helped us to clear up some
confusions in our understanding of the works of Markl-Papadima and Papadima.

\section{Background}

\subsection{The works of Markl-Papadima and Papadima}

We first summarize briefly the works of Markl-Papadima \cite{MaPa} and 
Papadima \cite{Pa}. See also \cite{BePa}. 

Let $X$ and $Y$ be the $\Bbb Q$-vector spaces spanned by $\{x_1,\dots,x_n\}$
and $\{y_1,\dots,y_n\}$, respectively. Let $\Bbb L_{\ast}^{\ast}=\Bbb L(X
\oplus Y)$ be the free bigraded $\Bbb Q$-Lie algebra generated by $X\oplus Y$, where  
the upper degree is given by the bracket length and the lower degree comes 
from setting $|X|=0$ and $|Y|=1$. Let $\hat{\Bbb L}=\hat{\Bbb L}_{\ast}$ be 
the completion 
of $\Bbb L_{\ast}^{\ast}$ with respect to the upper degree, i.e. with respect 
to the filtration $F_k{\Bbb L}_\ast={\Bbb L}_{\ast}^{\geq k}$. We have an 
induced filtration $F_k\hat{\Bbb L}$ for $\hat{\Bbb L}$.

Consider $\hat M=\text{Der}_{-1}^+\hat{\Bbb L}$, the vector space of
continuous derivations $\partial$ of $\hat{\Bbb L}$, which are homogeneous of 
lower degree $-1$ and have the property that $\partial F_k\hat{\Bbb L}\subset
F_{k+1}\hat{\Bbb L}$, $\forall k$, filtrated by 
$$F_s\hat M=\{\partial\,\,|\,\,\partial F_k\hat{\Bbb L}\subset F_{k+s}\hat{\Bbb L},
\forall k\}.$$

Since $\hat{\Bbb L}$ is free, a derivation
$\partial\in \text{Der}_{-1}^+\hat{\Bbb L}$ is completely determined by its
restriction to the free generators. Together with the assumption that 
$\partial$ is of lower degree $-1$ and $\partial F_k\hat{\Bbb L}\subset
F_{k+1}\hat{\Bbb L}$, $\forall k$, we may identify $\partial$ with a series
$$\partial=\partial_1+\partial_2+\cdots+\partial_s+\cdots,\quad \partial_s
\in\text{Hom}
(Y,{\Bbb L}^{s+1}_X)),\,\, \forall s
$$
where $\hat{\Bbb L}_X$ is the completed free $\Bbb Q$-Lie algebra on $X$, 
grade by the
bracket length (upper degree) and filtrated by $F_s\hat{\Bbb L}_X=
\hat{\Bbb L}_X^
{\geq s}$. We have
$$
\partial\in F_s\hat M\quad\Longleftrightarrow\quad\partial_{< s}=0,
$$
and obviously,
\begin{equation}
\hat M/F_s \hat M \cong \text{Hom}(Y,\hat{\Bbb L}_X^{\geq 2}/\hat{\Bbb L}_X^
{\geq s+1}).
\end{equation}

Consider the group $U$ of continuous lower degree zero Lie algebra 
automorphisms $u$ of $\hat{\Bbb L}$ with the property that $\text{gr}^1u=\text
{id}$ (with respect to the filtration $F_s\hat{\Bbb L}$. It acts on $\hat M$ by
conjugation and the action preserves the filtration $F_s\hat M$. Define
$${\cal M}_s=U \backslash\hat M\slash F_s\hat M$$
be the double quotient of $\hat M$ by $F_s\hat M$ and the $U$-action. We have
a natural projection
$$\text{pr}_s^{s+1}:{\cal M}_{s+1}\longrightarrow{\cal M}_s.$$

\medskip
\noindent{\bf Remark.} It is very important that the space ${\cal M}_s$ is 
the orbit space of the
vector space $\hat M\slash F_s\hat M$ under the action of $U$ and {\it not} 
the quotient vector space. In fact, the latter is trivial. Nevertheless, some 
linear structures do survive in ${\cal M}_s$. For example, there
is a natural scalar
multiplication on ${\cal M}_s$. For each $s$, there is 
a unique element $0\in{\cal M}_s$ invariant under the scalar multiplication 
and preserved by the projection 
$\text{pr}_s^{s+1}$. Moreover, $(\text{pr}_s^{s+1})^{-1}(0)\subset
{\cal M}_{s+1}$
is a vector space for each $s$. See the discussion at the beginning of 
Section 5.
\medskip
       
Now consider a classical link $L\in S^3$, consisting of $n$ disjointly embedded
circles, which are ordered and oriented. 
Let $G=\pi_1(S^3\setminus L)$ be the fundamental group of the link complement
of $L$, and $\Gamma_s G$ be the lower central series of $G$ with $\Gamma_1G=G$
and $\Gamma_s=(\Gamma_{s-1},G)$, where $(\,,\,)$ stands for the group 
commutator. The existence of a system of meridians and 
longitudes 
$m_i,l_i\in G$ ($1\leq i\leq n$) unique modulo conjugations $g_im_ig_i^{-1},
g_il_ig_i^{-1}$ by $g_i\in G$ ($1\leq i\leq n$) determines a {\it peripheral
structure} on $G$. The tower of nilpotent groups $G/\Gamma_{s+1} G$ inherits
such a peripheral structure. The structured tower of nilpotent groups 
$G/\Gamma_{s+1} G$ turns out to be an $I$-equivalence invariant of links
(\cite{St}). Here two links $L$ and $L'$ are {\it $I$-equivalent} if there 
is a 
proper (not necessarily locally flat) imbeddings of 
$S^1\times I_1\coprod\cdots\coprod S^1\times I_n$, $I_i=
I=[0,1]$, into $S^3\times I$ such that $L=S^1\times\{0\}\coprod\cdots\coprod
S^1\times\{0\}\subset S^3\times\{0\}$ and $L'=S^1\times\{1\}\coprod\cdots
\coprod
S^1\times\{1\}\subset S^3\times\{1\}$. In particular, 
every knot is $I$-equivalent to the unknot. 

The nilpotent group $G/\Gamma_{s+1} G$ admits a very simple presentation. Let $
{\Bbb F}={\Bbb F}(x_1,\dots,x_n)$ be the free group generated by 
$x_1,\dots,x_n$. 

\begin{theorem} {\rm (Milnor \cite{Mi2})} Given an $n$-component link $L$ 
with the link group $G=\pi_1
(S^3\setminus L)$. For every $s\geq1$, there exist elements $l_i^{(s)}\in{\Bbb
F}$ ($1\leq i\leq n$), such that
\begin{itemize}
\item $l_i^{(s+1)}\equiv l_i^{(s)}$ mod $\Gamma_s{\Bbb F}$, $1\leq i\leq n$;
\item $G/\Gamma_{s+1} G$ is isomorphic to 
\begin{equation}
\{\, x_1\dots,x_n\,\,|\,\,(l_i^{(s)},x_i),1\leq i\leq n,\,\Gamma_{s+1}
{\Bbb F}\,\,\}
\end{equation}
as structured groups where the peripheral structure of the latter group is 
given by $x_i,l_i^{(s)}$, $1\leq i\leq n$.
\end{itemize}
\end{theorem}      

The {\it Campbell-Hausdorff representation} is an injective homomorphism
$$\rho:{\Bbb F}(x_1,\dots,x_n)\longrightarrow\hat{\Bbb L}_X$$
into the Campbell-Hausdorff group of $\hat{\Bbb L}_X$ given by 
$\rho(x_i)=x_i$, $1\leq i\leq n$ and the Campbell-Hausdorff formula
$$\rho(gh)=\rho(g)+\rho(h)+\frac12[\rho(g),\rho(h)]+\cdots,\quad g,h
\in{\Bbb F}$$
where $[\,,\,]$ is the algebraic commutator (see \cite{Se} for an explicit 
form of the Campbell-Hausdorff formula).
For an $n$-component link $L$ with the 
structured group $G$, use the presentation (2.2), one can define a homomorphism
$\partial^{(s)}\in\text{Hom}(Y,\hat{\Bbb L}^{\geq 2}_X)$ by
$$\partial^{(s)}(y_i)=[\rho(l^{(s)}_i),x_i],\quad 1\leq i\leq n.$$
Under the natural identification
(2.1), we have a projection
$$\text{Hom}(Y,\hat{\Bbb L}^{\geq 2}_X)
\longrightarrow\text{Hom}(Y,\hat{\Bbb L}^{\geq2}_X/\hat
{\Bbb L}^{\geq s+1}_X)\longrightarrow{\cal M}_s.$$
Then, define $CH^{(s)}(L)\in{\cal M}_s$ to be the image of $\partial^{(s)}$
under this projection. We have

\begin{theorem}{\rm (Papadima \cite{Pa})} $CH^{(s)}(L)\in{\cal M}_s$ depends 
only on the $I$-equivalence class of $L$. Moreover, we have
$$\text{\rm pr}_s^{s+1}(CH^{(s+1)}(L))=CH^{(s)}(L).$$
\end{theorem}

\subsection{Finite type link invariants} 

Here we briefly review the results in \cite{Li}. 

First, we will refer the reader to the standard references \cite{Ba1,Bi,BiLi} 
for the definition finite type knot invariants (or Vassiliev invariants). The 
definition can be easily generalized to links and string links \cite{Ba2,Li}.

Briefly, an {\it $n$-component string link $\sigma$} is 
a self-concordance of $n$ fixed points in the plane 
${\Bbb R}^2$. For each $s$, such a string link $\sigma$ induces an 
automorphism
$$\phi_\sigma^{(s)}:{\Bbb F}/\Gamma_{s+1}{\Bbb F}\longrightarrow
              {\Bbb F}/\Gamma_{s+1}{\Bbb F}$$
given by
$$\phi^{(s)}_\sigma(x_i)=l_i^{(s)}x_i(l_i^{(s)})^{-1},\quad 1\leq i\leq n$$
for certain $l_i^{(s)}\in{\Bbb F}$ and
$$l_i^{(s+1)}\equiv l_i^{(s)}\quad\text{mod}\quad \Gamma_s{\Bbb F}.$$
 
Let $\hat{\Bbb P}=\hat{\Bbb P}(X_1,\dots,X_n)$ be the ring of formal 
power series (with rational coefficients)
in non-commutative variables $X_1,\dots,X_n$. It comes naturally with a 
grading and a filtration
$F_s\hat{\Bbb P}$. An {\it expansion} of the free
group ${\Bbb F}={\Bbb F}(x_1,\dots,x_n)$ is an injective homomorphism $E:{\Bbb
F}\rightarrow\hat{\Bbb P}$ such that
$$E(x_i)=1+X_i+\text{terms of order $\geq 2$}.$$

Let ${\Bbb P}^{(s)}=\hat{\Bbb P}/F_{s+1}\hat{\Bbb P}$. The automorphism 
$\phi_\sigma^{(s)}$ induced by a string link $\sigma$ can be converted into an
automorphism
$$\Phi_\sigma^{(s)}:{\Bbb P}^{(s)}\longrightarrow{\Bbb P}^{(s)}$$
using the expansion $E$.

Suppose now we have a singular string link $\tau$ with $k\geq s$ double 
points. 
Let $\sigma_j$, $j=1,2,\dots,2^k$  
be the $2^k$ string links obtained by resolving the $k$ double points on 
$\tau$ into positive or negative crossings. Let $\epsilon_j=1$ or $-1$ 
depending on whether $\sigma_j$ has even or odd number of negative 
resolutions. Then, we say that the string link invariant $\Phi_\sigma^{(s)}$
is of finite type in the following sense.
 
\begin{theorem}{\rm (Lin \cite{Li})} In $\text{\rm End}({\Bbb P}^{(s)})$, 
we have
$$\sum_{j=1}^{2^k}\epsilon_j\Phi_{\sigma_j}^{(s)}=0.$$ 
\end{theorem} 

\section{Campbell-Hausdorff invariants are of finite type}

We now use the {\it Campbell-Hausdorff expansion} $E:{\Bbb F}\rightarrow
\hat{\Bbb P}$ given by
$$E_{CH}(x_i)=\text{exp}(X_i)=\sum_{s=0}^\infty\frac1{s!}X_i^s.$$
For a string link $\sigma$, the induced ring automorphism $\Phi_\sigma^{(s)}$
via $E_{CH}$ is completely determined by $\Phi_\sigma^{(s)}(X_i)$, 
$1\leq i\leq n$, which are Lie elements in ${\Bbb P}^{(s)}$. We may relate
these Lie elements with $[\rho(l_i^{(s)}),x_i]$, $1\leq i\leq n$ in the 
following way.

We have 
\begin{equation}
\Phi_\sigma^{(s)}(X_i)\equiv\log\left[E_{CH}(l_i^{(s)})E_{CH}(x_i)
E_{CH}((l_i^{(s)})^{-1})\right]\quad\text{mod}\,\, F_{s+1}\hat{\Bbb P}.
\end{equation}
Obviously, 
\begin{equation}
\log E_{CH}((l_i^{(s)})^{\pm1})=\pm\rho(l_i^{(s)})
\end{equation} 
with each $x_j$ in the right-hand side being replaced by $X_j$.
From now on, we will abuse the notation by identifying these two Lie elements. 
In general, we will not distinguish Lie elements in $x_j$'s or in $X_j$'s.

\begin{lemma} With the same notation as in Theorem 2.3, let $
l_i^{(s)}(\sigma_i)$ be the $i$-th longitude of the string link $\sigma_j$.
Then
$$\sum_{j=1}^{2^k}\epsilon_j[\rho(l_i^{(s)}(\sigma_j)),X_i]
\equiv 0\quad\text{\rm mod}\,\,F_{s+1}\hat{\Bbb L}_X.$$
\end{lemma}

\begin{proof} We argue first inductively that there is a basis for $\hat
{\Bbb L}_X
/F_{s+1}\hat{\Bbb L}_X$, such that for a given base element $e$, 
its coefficient 
in $\rho(l_i^{(s)})$ can be expressed as a polynomial of coefficients of base 
elements or degree $\leq\text{degree}(e)+1$ in $\Phi_\sigma^{(s)}(X_i)$. 

Suppose we have chosen a basis for $\hat{\Bbb L}_X
/F_{s+1}\hat{\Bbb L}_X$ up to degree $r+1$ such that for a given 
base element $e$ of degree $\leq r$, its coefficient in $\rho(l_i^{(s)})$ is a
polynomial of coefficients of base 
elements or degree $\leq\text{degree}(e)+1$ in $\Phi_\sigma^{(s)}(X_i)$.
We want to extend this basis to a basis up to degree $r+2$ and determine the
coefficients of base elements of degree $r+1$ in $\rho(l_i^{(s)})$. 
   
Notice the following two things:
\begin{enumerate}
\item for $e\in\hat{\Bbb L}_X$, $[e,X_i]=0$ iff $e$ is a multiple of $X_i$; and 
\item by the choice of longitude, $\rho(l_i^{(s)})$ contains no $X_i$.
\end{enumerate}

By 1., we may pick up degree $r+1$ base elements $e$'s not equal to $X_i$ and
to extend $[e,X_i]$'s to a basis of degree $r+2$. By 2., the coefficients of 
base elements of degree $r+1$ in $\rho(l_i^{(s)})$ are all lifted up to become
coefficients of base elements of degree $r+2$. Then we may use (3.1), (3.2) 
and the Campbell-Hausdorff formula to finish the induction.

Now we see that there is a basis for $\hat{\Bbb L}_X
/F_{s+1}\hat{\Bbb L}_X$ such that the coefficient of a given base element
in $[\rho(l_i^{(s)}),X_i]$ is a polynomial of coefficients of base elements
in $\Phi_\sigma^{(s)}(X_i)$. The lemma then follows from Theorem 2.3 
and the fact that 
polynomials of finite type invariants are still of finite type 
(with compatible order -- the maximal number of double points singular links 
can have so that the invariant does not necessarily vanish on derived 
alternating sum of links).
\end{proof}  

We come to our main result of this section now. It says that the 
Campbell-Hausdorff invariants are of finite type with respect to 
marked singular links.

\begin{definition}  A marked singular link is a singular link in the usual 
sense with a marked point on each component away from double points. 
An isotopy of a marked singular 
link is an usual isotopy with the restriction that the marked points are
not allowed to pass through double points. 
\end{definition}

Let $J$ be a marked singular link with $k$ double points. We may find a 
disk $D$ in the $S^3$ which
intersects $J$ transversally at the marked points. This converts $J$ into a 
singular string link. By picking a base point on $D$ and paths on $D$ 
connecting the base point and the marked points on $J$, we also get a 
{\it consistent} choice of 
meridian-longitude pairs for all links $J_1, 
J_2,\dots,J_{2^k}$ derived from $J$ by resolving  
the double points. Such a choice of meridian-longitude pairs is consistent
in the sense that if we choose another disk $D'$, according to the basic
construction of \cite{Pa}, there will be an element $u\in U$ which sends
derivations associated with all $J_j$ determined by $D$ to that 
determined by $D'$. This will certainly not be true if marked points are
allowed to pass through double points. 

Thus, we get a well defined invariant of marked singular links 
$$CH^{(s)}(J)\in{\cal M}_s$$
which is the image of
$$\sum_{j=1}^{2^k}\epsilon_j\partial^{(s)}(J_j)\in\text{Hom}(Y,
\hat{\Bbb L}_X^{\geq2}),$$      
determined by the consistent meridian-longitude pairs coming from $D$,
under the projection
$$\text{Hom}(Y,
\hat{\Bbb L}_X^{\geq2})\longrightarrow\text{Hom}(Y,\hat{\Bbb L}_X^{\geq 2}/
\hat{\Bbb L}_X^{\geq s+1})\longrightarrow{\cal M}_s.$$

\begin{theorem} Let $J$ be a marked singular link with $k\geq s$ double 
points. Then
$$CH^{(s)}(J)=0.$$
\end{theorem}

\begin{proof} We realize $J$ as the closure of a singular string link $\tau$
with $k$ double points. Then $L_j$'s are closures of resolutions of $\tau$.
For a fixed consistent choice of longitudes of $L_j$'s, Lemma 3.1 says that 
the alternating sum of
the resulting derivations is zero in $\hat M/F_s\hat M$. 
Therefore, $CH^{(s)}(J)=0$.
\end{proof}

Thus, the invariant of marked singular links derived from
$CH^{(s)}$ vanishes for every marked singular links with more than 
$s-1$ double points. It is in this sense we say that the Campbell-Hausdorff 
invariants are of {\it finite type with respect to marked singular links}.
See some further discussion at the beginning of Section 5.

\section{Quantum invariants can not detect non-invertibility}

A link $L$ is {\it non-invertible} if it is not isotopic to the link $L'$ 
obtained by reversing the orientation of each component of $L$. It is a 
well-known fact that quantum invariants can not detect the 
non-invertibility of 
a link. We include here a proof of this fact since we could not find one
in the literature. The basic ingredients are
all in \cite{ReTu1, ReTu2} already and the reader is referred to there for
detailed discussions of the construction of quantum invariants.

Let $A=U_q\frak g$
be the quantized universal enveloping algebra of a simple Lie algebra 
$\frak g$. Together with an invertible element $R\in A\otimes A$ and an
invertible central
element $v\in A$, $(A,R,v)$ is a ribbon Hopf algebra. There is also 
an invertible element $u\in A$ which, among other things, has the following 
property useful to us. 

Let $(V,\rho)$ be an $A$-module and $V^{\lor}$ be its dual. If we identify
$V^{\lor\lor}$ with $V$ via the canonical isomorphism, we have 
$$\rho^{\lor\lor}(a)=\rho(u)\rho(a)\rho(u^{-1}), \quad \forall a\in A.$$

We will denote by
$V^{\lambda}$ the finite dimensional irreducible representation of $U_h\frak g$
with the highest weight $\lambda$. The key fact we need is that there
is an $A$-linear isomorphism
$$w_{\lambda}:(V^{\lambda})^{\lor}\longrightarrow V^{\lambda^*}$$ 
such that $(\lambda^*)^*=\lambda$ and
$$w_\lambda^{\lor}\cdot (w_{\lambda^*})^{-1}:V^{\lambda}\longrightarrow
(V^\lambda)^{\lor\lor}=V^\lambda$$
is the multiplication by $uv^{-1}$.

Now let $B$ be a braid with $m$ strands. Suppose that all the strands of $B$
point downward. We denote by $B^-$ the braid obtained from $B$ by reversing
the arrow on each strand of $B$. Then we color each strand of $B$ by an
$A$-module $V^\lambda$, in such a 
way that the functor $F$ defined in Theorem 5.1 of \cite{ReTu1} applied to 
the colored braid $B_{\lambda_1,\dots,\lambda_m}$ gives an automorphism
$$F(B_{\lambda_1,\dots,\lambda_m}):V^{\lambda_1}\otimes\cdots\otimes V^{\lambda_m}
\longrightarrow V^{\lambda_1}\otimes\cdots\otimes V^{\lambda_m}.$$
By Lemma 5.1 of \cite{ReTu2}, for
$$F(B^-_{\lambda_1^*,\dots,\lambda_m^*}):
(V^{\lambda_1^*})^\lor\otimes\cdots(V^{\lambda^*_m})^\lor
\longrightarrow(V^{\lambda_1^*})^\lor\otimes\cdots(V^{\lambda^*_m})^\lor,
$$ 
we have
$$
F(B^{-}_{\lambda^*_1,\dots,\lambda^*_m})=  
w_{\lambda^*_1}^{-1}\otimes\cdots w_{\lambda^*_m}^{-1}\cdot
F(B_{\lambda_1,\dots,\lambda_m})\cdot w_{\lambda^*_1}\otimes\cdots\otimes
w_{\lambda^*_m}.$$

Conjugate further by $w_{\lambda_1}^\lor\otimes\cdots\otimes w_{\lambda_m}^\lor
$, we get
$$F(B^-_{\lambda_1^\lor,\dots,\lambda_m^\lor})=uv^{-1}\otimes\cdots\otimes
uv^{-1}\cdot F(B_{\lambda_1,\dots,\lambda_m})\cdot
vu^{-1}\otimes\cdots\otimes vu^{-1},$$
where 
$$F(B^-_{\lambda_1^\lor,\dots,\lambda_m^\lor})=
F(B^-_{\lambda_1,\dots,\lambda_m})^{\lor}.$$

Now let $\text{Tr}_q$ be the quantum trace, i.e.
$$\text{Tr}_q(a)=\text{Tr}\,(uv^{-1}a), \quad \forall a\in A$$
where $\text{Tr}\,(\cdot)$ is the ordinary trace. Then, we have
$$
\begin{aligned}
\text{Tr}_q(F(B_{\lambda_1,\dots,\lambda_m}))
  &=\text{Tr}_q(F(B^-_{\lambda_1,\dots,\lambda_m})^\lor) \\
  &=\text{Tr}\,((uv^{-1})^\lor F(B^-_{\lambda_1,\dots,\lambda_m})^\lor) \\
  &=\text{Tr}\,(uv^{-1}F(B^-_{\lambda_1,\dots,\lambda_m})) \\
  &=\text{Tr}_q(F(B^-_{\lambda_1,\dots,\lambda_m})).
\end{aligned}
$$

The value of the quantum invariant on the colored link 
obtained as the closure of
$B_{\lambda_1,\dots,\lambda_m}$ is simply $\text{Tr}_q
(F(B_{\lambda_1,\dots,\lambda_m}))$ together with a normalization factor 
independent of the overall orientation.
Therefore, such a quantum invariant can not detect the non-invertibility of 
a link since it has the same value on the closures of $B_{\lambda_1,\dots,
\lambda_m}$ and $B^-_{\lambda_1,\dots,\lambda_m}$.
  
\section{An example and some comments}

Since the action of $U$ on $\hat M/F_s\hat M$ is unipotent, its restriction on
$F_{s-1}\hat M/F_s\hat M$ is trivial. Thus, just like in the
case of $\bar\mu$-invariant, there is no indeterminacy for the 
Campbell-Hausdorff invariants when 
$$\rho(l_i^{(s)})\in\hat{\Bbb L}_X^{\geq s-1},\quad 1\leq i\leq n.$$ 
In other words, we have an injection
\begin{equation}\text{Hom}(Y,\hat{\Bbb L}_X^{\geq s}/\hat{\Bbb L}_X^{\geq s+1})
\hookrightarrow{\cal M}_s.\end{equation}
Actually, when $J$ is a marked singular link with $s-1$ double points, it is
easy to see that 
$$CH^{(s)}(J)\in\text{Hom}(Y,\hat{\Bbb L}_X^{\geq s}/\hat{\Bbb L}_X^
{\geq s+1}).
$$ 
According to Theorem 3.1, this gives us a weight system on
{\it marked chord diagrams} with $s-1$ chords. Here a marked chord diagram
is a usual chord diagram on oriented circles with an additional mark on each 
circle away from the end points of the chords.   
Compare with the conventional 
theory of finite type invariants \cite{Ba1}, it is interesting to see that
such a weight system can be integrated into a link invariant not living in
a vector space. It is certainly desirable to explore such a phenomenon
more systematically.

Like $\bar\mu$-invariants again, (5.1) makes the computation of 
non-trivial Campbell-Hausdorff invariants of lowest order quite easy once we 
have the 
presentation (2.2). In fact, at this level, the Campbell-Hausdorff invariants
are exactly the same as $\bar\mu$-invariants. So we may use the existing 
computation of $\bar\mu$-invariants in the literature (see \cite{Co}, for 
example). 

Let us now try to understand how the longitudes change when we reverse the
orientation of a link. 

\begin{lemma} Let $L'$ be obtained from $L$ by reversing the orientation. 
We may have a system of longitudes $\l_i^{(s)}\in{\Bbb F}$, $1\leq
i\leq n$, for $L$ such that $\bar l_i^{(s)}\in{\Bbb F}$, $1\leq i\leq n$, which
are obtained by reading $l_i^{(s)}$ backward, $1\leq i\leq n$ respectively, 
as a word in $x_1,\dots,x_n$, form a system of longitudes
for $L'$. 
\end{lemma}

\begin{proof} We may use the Wirtinger presentation of the link group to 
get longitudes $l_i^{(s)}\in{\Bbb F}$. See \cite{Mi2}. Simply note that when 
we change $L$ to $L'$ 
by reverse the orientation, a positive (negative) crossing in a projection of 
$L$ remains 
to be positive (negative) in $L'$. Therefore, when we read the geometric 
longitudes of $L$ backward as words in Wirtinger generators, we get the 
geometric longitudes of $L'$. Now the process of reducing geometric longitudes
of $L$ to $l_i^{(s)}$ can also be read backward to get the reduction of 
geometric longitudes of $L'$ to $\bar l_i^{(s)}$. This finishes the proof.
\end{proof}

\medskip

\noindent{\bf Example 5.1.} Let $L$ be the Borromean rings with appropriate
orientation and numbering of components such that
$$\bar\mu(123)=\bar\mu(231)=\bar\mu(312)=-\bar\mu(213)=-\bar\mu(321)=-\bar\mu
(132)=1.$$
This implies $CH^{(3)}(L)=-CH^{(3)}(L')\neq 0\in {\cal M}_3$, or 
$$CH^{(3)}(L)\neq CH^{(3)}(L').$$
Thus, the link invariant $CH^{(3)}$ can detect the
non-invertibility of a link and it can not be determined by quantum invariants.

\medskip

The non-invertibility of a knot is much subtler to detect. 
The first non-invertible knot was detected by Trotter in 1963 \cite{Tr}.
It is also known that there are knots non-invertible up to concordance
\cite{Liv}. More recently, Kuperberg \cite{Ku} showed that the Conway's 
11-crossing
knot is non-invertible by counting epimorphisms of the knot group of
Conway's knot into the sporadic simple group $M_{11}$ of order 7920. Hyperbolic
structures on knot complements may also be used to detect non-invertibility
of knots \cite{Ka}.

For a knot $K$ with the group $G=\pi_1(S^3\setminus K)$, simple homological
reasons (Stallings' theorem \cite{St}) imply
$$\Gamma_2G=\Gamma_3G=\cdots=\Gamma_s G=\cdots.$$
In particular, the longitude of a knot can not be seen in the quotients of $G$
by its lower central series. For that reason and the nature of Vassiliev 
knot invariants as well as the current work, we feel more convinced now
of the truth of the following conjecture which we made several years ago 
when people 
started to talk about Vassiliev invariants. 

\begin{conjecture} Vassiliev knot invariants can not detect the 
non-invertibility of a knot.
\end{conjecture}

\end{document}